\documentclass[12pt]{amsart}
\usepackage{amssymb,amsfonts}
 \usepackage{a4wide}
\usepackage{graphicx}
\usepackage{psfig}
\usepackage{amsrefs}

\numberwithin{equation}{section}

\newtheorem{theorem}{Theorem}

\newtheorem{thm}{Theorem}[section]
\newtheorem{lemma}[thm]{Lemma}

\theoremstyle{definition}
\newtheorem{defn}{Definition}

\newtheorem{rmk}[thm]{Remark}

\newcommand{\R}{{\mathbb R}}

\newcommand{\N}{{\mathbb N}}

\renewcommand{\hat}{\widehat}
\renewcommand{\tilde}{\widetilde}

\newcommand{\sgn}{\operatorname{sgn}}

\begin{document}

\title{The Lorenz attractor is mixing}

\author{Stefano Luzzatto}
\address{Mathematics Department, Imperial College, 180 Queen's Gate,
London SW7 2AZ}
\email{stefano.luzzatto@imperial.ac.uk}
\urladdr{www.ma.ic.ac.uk/~luzzatto}
\author{Ian Melbourne}
\address{
Department of Mathematics and Statistics,  Ê 
 School of ECM, Ê 
 University of SurreyÊ 
 Guildford, SurreyÊ 
 GU2 7XH, UK}
 \email{ism@login.math.uh.edu}
     \urladdr{www.maths.surrey.ac.uk/showstaff?I.Melbourne}
\author{Frederic Paccaut}
\address{LAMFA UMR 6140 du CNRS, 
University of Picardie, 
 33, rue Saint-Leu
 80039 Amiens cedex 1, France}
\email{frederic.paccaut@u-picardie.fr}
\urladdr{lamfa.u-picardie.fr/paccaut/}

\begin{abstract}
We study a class of geometric Lorenz flows, introduced independently
by Afra{\u\i}movi\v{c}, Bykov \& Sil{\cprime}nikov and by Guckenheimer
\& Williams, and give a verifiable condition for such flows to be mixing.
As a consequence, we show that the classical Lorenz attractor is mixing.
 \end{abstract}

\maketitle

\section{Introduction and statement of results}
\subsection{The Lorenz equations}

Many systems of nonlinear differential equations 
that were first studied almost 50 years ago and which were
motivated mainly by problems in 
geophysical and astrophysical fluid dynamics and dynamical 
meteorology \cites{Bul55, Rik58,All62, Hid04},  
remain difficult to understand rigorously to the present day.
In 1963, Lorenz~\cite{Lor63} introduced the following
system of differential equations:
\begin{align} \label{eq-lorenz}  \nonumber
\dot{x} &= 10 (y-x)  \\
\dot{y} &=  28 x -y -xz     \\
\dot{z} &= xy - {\textstyle\frac{8}{3}} z  \nonumber
\end{align}
Approximate numerical studies of these equations led Lorenz to 
emphasise the possibility and importance of
\emph{sensitive dependence on initial conditions} even in such 
 simplified models of natural phenomena. 
A combination of results obtained 
over the last 25 years~\cites{Afr77, Wil79, GucWil79, Pes92, Sat92}
and culminating in the work of Tucker~\cites{Tuc99,Tuc02} gives the following 
statement (see \cites{Spa82, Via01} for detailed surveys): 

\begin{quote}
    \emph{The Lorenz equations admit a robust attractor \( 
    A \) which supports a ``physical'' ergodic invariant probability 
    measure \( \nu \) with a positive Lyapunov exponent}
\end{quote}  
Recall that the measure $\nu$
is called physical, or Sinai-Bowen-Ruelle (SRB), if 
for Lebesgue almost every solution $u(t)\in\R^3$ starting close to $A$ and all 
continuous functions $h:\R^3\to\R$,
\[
\lim_{T\to\infty}\frac1T \int_0^T h(u(t))\,dt = \int_A h\,d\nu.
\]
In this paper we  take a further step in the understanding of the 
statistical properties of the Lorenz attractor. 
A measure \( \nu \) is mixing for a flow $\Phi_t$ if 
\[ 
\nu(\Phi_{t}(A)\cap B) \to \nu(A) \nu(B)
\]
for all measurable sets \( A, B  \), as \( 
t\to \infty \). We say that the Lorenz attractor is mixing if the 
SRB measure \( \nu \) mentioned above is mixing.
We prove the following 

\begin{theorem}\label{thm-lorenz}
    The Lorenz attractor is mixing.
\end{theorem}

In fact, our result shows that the Lorenz attractor is \emph{stably mixing}:
sufficiently small $C^1$ perturbations of the flow are mixing. 

There are relatively few explicit examples of flows that have been proved to be
mixing.    Anosov~\cite{Ano67} showed that geodesic flows for compact manifolds
of negative curvature are mixing, and this was generalised~\cite{BurKat94} 
to include contact flows.  Moreover, mixing persists under $C^1$ perturbations.
For codimension one Anosov flows~\cite{Pla72} and for Anosov flows with a 
global infranil cross-section~\cite{Bow76}, the set of mixing flows is 
$C^1$ open, but the corresponding result for general Anosov flows is not known.
Recently~\cite{FieMelTor04} it was shown that mixing holds for a $C^1$-open and 
$C^r$-dense set of $C^r$ Anosov flows for all $r\ge1$.
(In these references, mixing is proved for any equilibrium measure
for a H\"older potential.)
However, the conditions for stable mixing in~\cites{FieMelTor03,FieMelTor04} for Anosov and 
uniformly hyperbolic (Axiom~A flows) are not easily verifiable.

The Lorenz attractor is an example of a {\em singular hyperbolic attractor}~\cite{MorPacPuj99}
(uniformly hyperbolic, except for a singularity due to the attractor containing
an equilibrium).
Somewhat surprisingly, we show that the singular nature of the Lorenz attractor
assists in the search for a verifiable condition for mixing.

\subsection{Geometric Lorenz attractors}

Afra{\u\i}movi\v{c}, Bykov and Sil{\cprime}nikov~\cite{Afr77} and Guckenheimer
and Williams~\cites{GucWil79,Wil79} introduced a geometric
model that is an abstraction of the 
numerically-observed features possessed by solutions to~\eqref{eq-lorenz}.
Tucker~\cites{Tuc99,Tuc02} proved that the geometric model is valid,
so the Lorenz equations define a geometric Lorenz flow.

Accordingly, our approach to Theorem~\ref{thm-lorenz} is to establish mixing for
geometric Lorenz flows satisfying certain hypotheses, and then to verify
from~\cites{Tuc99,Tuc02} that the hypotheses are satisfied for the
Lorenz equations.  Roughly speaking, a {\em geometric Lorenz flow} is 
the natural extension of a {\em geometric Lorenz semiflow}
which is itself a suspension flow built over a certain type of
one-dimensional expanding map $f$.  Precise definitions
are given in Section~\ref{sec-tucker}.   In the literature, a standard 
assumption is that the map $f$ is {\em locally eventually onto
(l.e.o.)}, see Definition~\ref{def-leo}.  
We prove that this is a sufficient condition for the 
corresponding flow to be mixing.

\begin{theorem}\label{thm-geom}
Let \( f_t \) be a geometric Lorenz flow.    Suppose that
the associated one-dimensional map $f$ is l.e.o.
Then \( f_t \)  is mixing (and even Bernoulli).
\end{theorem}

The remainder of the paper is structured as follows.
In Section~\ref{sec-tucker}, we discuss the geometric Lorenz flow
and its relations to the Lorenz attractor.
In particular, we indicate how Theorem~\ref{thm-lorenz} follows from 
Theorem~\ref{thm-geom}.
In Section~\ref{sec-mix}, we prove Theorem~\ref{thm-geom}.
In Section~\ref{sec-ext}, we discuss extensions of our
main results and some related future directions.

 \section{The Lorenz attractor is a geometric Lorenz flow}
 \label{sec-tucker}
 
 In this section we collect and organise several results from the 
 existing literature on the relation between the Lorenz attractor 
 and geometric Lorenz flows. 
  We recall first of all some basic relevant facts about the Lorenz 
equations~\eqref{eq-lorenz}, see \cite{Spa82}.  The origin is an 
 equilibrium of saddle type with two negative (stable) and one 
 positive (unstable) eigenvalues $\lambda_{ss}<\lambda_s<0<\lambda_u$.
It is also the case that $\lambda_u>|\lambda_s|$.

Suppose that a finite number of nonresonance conditions are satisfied
so that the vector field is smoothly linearisable in a neighbourhood
of $0$.   In these coordinates, the flow near $0$ has the form
$(x_1,x_2,x_3)\mapsto
(e^{\lambda_u t}x_1, e^{\lambda_{ss} t}x_2, e^{\lambda_s t}x_3)$.
By a linear rescaling, we can suppose that the domain of
linearisation of the flow includes the cube $[-1,1]^3$.
Define $\Sigma=\{(x_1,x_2,x_3): |x_1|,|x_2|\le1,\,x_3=1\}$ and
$\tilde\Sigma=\{(x_1,x_2,x_3): x_1=\pm1,\,|x_2|,|x_3|\le1\}$.
Then we define the first hit map $P_0:\Sigma\to\tilde\Sigma$ by
\[
P_0(x_1,x_2,1)=(e^{\lambda_u r_0}x_1,e^{\lambda_{ss} r_0}x_2,e^{\lambda_s r_0})=(\sgn x_1, \tilde x_2,\tilde x_3),
\]
where $r_0$ is the ``time of flight'' from $\Sigma$ to $\tilde\Sigma$.
Solving $e^{\lambda_u r_0}x_1=\sgn x_1$ for $r_0=-(\ln |x_1|)/\lambda_u$,
we obtain
\[
P_0(x_1,x_2,1)=(\sgn x_1,|x_1|^\beta x_2, |x_1|^\alpha),
\]
where $\alpha=|\lambda_s|/\lambda_u\in(0,1)$ and
$\beta=|\lambda_{ss}|/\lambda_u>0$.
Note that $P_0:\Sigma\to\tilde\Sigma$ is well-defined on
$\Sigma\setminus W^s(0)$ and that $r_0$ has a logarithmic singularity
at $x_1=0$.

Tucker \cites{Tuc99, Tuc02} proves that there exists a compact trapping region
$N\subset \Sigma$ such that
the \emph{Poincar\'e first return map} 
\(
P: N\setminus W^s(0) \to N
\)
is well defined.   We can decompose $P=P_1\circ P_0$ where $P_1$
is the first-hit map between $\tilde\Sigma$ and $\Sigma$.   Note that $P_1$
is a diffeomorphism where it is defined, and that the time of
flight $r_1$ for $P_1$ is bounded.
Hence the time of flight $r=r_0+r_1$ for the full return map $P$ is smooth
except for a logarithmic singularity at $x_1=0$.
Moreover, the following crucial \emph{hyperbolicity estimate} holds. 

\begin{lemma}[\cite{Tuc02}] \label{lem-cone}
The return map $P$ admits a forward invariant cone field.
In other words, there exists a cone $\mathcal{C}(u)$ inside $\Sigma$ at each
point $u\in N\setminus W^s(0)$ such that $(dP)_u\mathcal{C}(u)$ is strictly contained inside
$\mathcal{C}(Pu)$.

Moreover, there exist constants $c>0$, $\tau>1$, such that
for each $u\in N\setminus W^s(0)$, 
\[
\|(dP^n)_uv\|\ge c\tau^n\|v\|,
\]
for every $v\in \mathcal{C}(u)$ and $n\ge1$.
\end{lemma}

A consequence of Lemma~\ref{lem-cone} is that the return map $P$
has an invariant stable foliation with $C^{1+\varepsilon}$ leaves for some
$\varepsilon>0$.  
Let $I=[-1,1]$.
We obtain  a singular one-dimensional map $f:I\to I$
by quotienting along stable leaves.
At the same time, $r$ reduces to a singular map
$r:I\to\R^+$.
Let $J=I\setminus \{0\}$.  Then it is immediate that on $J$:
\begin{description}
 \parskip = 3pt 
\item[{\bf f}(i)] $f$ is $C^{1+\varepsilon}$. 
\item [{\bf f}(ii)]   $|f^{(n)}|\ge c\tau^n$ for all $n\ge1$.
\item[{\bf f}(iii)] $C^{-1}|x|^{\alpha-1}\le f'(x) \le C|x|^{\alpha-1}$,
\item[{\bf r}(i)] $r(x)\to\infty$ as $x\to0$, 
\item[{\bf r}(ii)] $|r(x)-r(y)|\le C|\ln |x|-\ln |y||$
for all $x,y>0$ and all $x,y<0$,
\end{description}
where $c>0$, $\tau>1$ are the constants in Lemma~\ref{lem-cone},
$\alpha=|\lambda_s|/\lambda_u\in(0,1)$, and $C\ge1$ is a constant.
Note that it follows from~{\bf r(ii)} that $r$ is Lebesgue integrable
and is $C^1$ on $J$.

We have described how to pass from the original three-dimensional flow
defined by the Lorenz equations~\eqref{eq-lorenz}
to a two-dimensional Poincar\'e map $P:\Sigma\to\Sigma$ and then
to a one-dimensional expanding map $f:I\to I$.
The process is reversible: taking the natural extension of $f$
recovers the Poincar\'e map $P$ and taking the suspension of the 
$P$ by the roof function $r$ recovers the original Lorenz flow.
Alternatively, these steps can be carried out
in the opposite order to recover the Lorenz flow as the natural extension
of the suspension semiflow of the map $f$ by the roof function $r$.
The latter viewpoint is the {\em geometric Lorenz flow} 
construction~\cites{GucWil79,Wil79}.
For completeness, the notions of suspension and natural extension are
recalled at the end of this section.

\begin{figure}[ht]
\centerline{%
\psfig{file=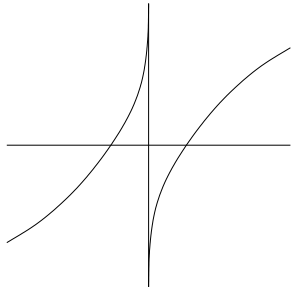,width=1.5in}
}
\caption{The graph of a Lorenz-like expanding map $f:[-1,1]\to[-1,1]$.}
\label{fig-exp}
\end{figure}

\begin{defn}  
Let $f:I\to I$.
Assume that $f(0)$ is undefined, with $f(0^+)=-1$ and $f(0^-)=+1$.
It is assumed that $f(1)\in(0,1)$ and $f(-1)\in(-1,0)$.
(See Figure~\ref{fig-exp}.)
If conditions~{\bf f(i)}--{\bf f(iii)} are satisfied, then $f$
is a {\em Lorenz-like expanding map}.
\end{defn}

\begin{defn}   A semiflow $f_t$ is called a {\em geometric Lorenz semiflow}
if it is the suspension over a Lorenz-like expanding map $f:I\to I$
by a roof function $r$ satisfying conditions~{\bf r(i)} and~{\bf r(ii)}.
A flow $f_t$ is called a {\em geometric Lorenz flow}
if it is the natural extension of a geometric Lorenz semiflow.
\end{defn}

The results of Tucker described above imply that the Lorenz equations
indeed define a geometric Lorenz flow.
Moreover, as outlined below, Tucker showed that the associated
expanding map $f:I\to I$
is transitive for the Lorenz equations, thus establishing the existence
of the Lorenz attractor.

In the following definition, we stress that $f(0)$ is undefined.
In particular, the endpoints $\pm1$ have no preimages.
\begin{defn}  
\label{def-leo}
The map $f:I\to I$ is {\em locally eventually onto
(l.e.o.)} if for any open set $U\subset J$, there exists $k\ge0$ such
that $f^kU$ contains $(0,1)$.
\end{defn}

\begin{rmk}   Since $f(0,1)$ contains $(-1,0)$ and vice versa,
we could equally use $(-1,0)$ instead of $(0,1)$
in the definition of l.e.o.
Clearly, the l.e.o.\ property implies that $f$ is topologically transitive.
The exact formulation of l.e.o.\ varies in the literature,
but our definition agrees with the one in~\cite{GleSpa93}.
\end{rmk}

\begin{rmk}
It was anticipated in~\cite{Wil79} that {\bf f(ii)} would hold
with $c=1$, $\tau>\sqrt 2$.    
It is easy to check that 
this is a sufficient condition for $f$ to be l.e.o.

Surprisingly~\cite{Tuc02}, for the actual Lorenz equations it turns out 
that $f'(x)<1$ for certain values of $x\in I$.
Nevertheless, it is still the case that $f$ is l.e.o., see Tucker~\cite{Tuc02}.
\end{rmk}

It can be shown that Lorenz-like expanding maps satisfying the l.e.o.\
condition have a unique ergodic probability measure $\mu$ that
is equivalent to Lebesgue (see for example Section~\ref{sec-mix}).
The suspended measure $\nu=\mu^r$ defines an SRB measure for the
geometric Lorenz flow.
In particular, Theorem~\ref{thm-lorenz} follows from Theorem~\ref{thm-geom}.
The latter result is proved in Section~\ref{sec-mix}.

\begin{rmk}   \label{rmk-Ratner}
The approach in this paper establishes weak mixing for
geometric Lorenz semiflows.   It then follows from
Ratner~\cite{Rat78} that such semiflows, and the corresponding flows,
are Bernoulli.  In particular, they are mixing.
\end{rmk}

\subsubsection*{Suspensions and natural extensions}
As promised, we end this section by recalling the notions of suspension
and natural extension.
Let $(X,\mu)$ be a probability space and
\mbox{$f:X\to X$} a (noninvertible) measure preserving transformation.
Let $r:X\to\R^+$ be an $L^1$ roof function.
Define the suspension $X^r=\{(x,u):x\in X,\,u\in[0,r(x)]\}/\sim$
where $(x,r(x))\sim (fx,0)$.
Form the suspension semiflow $f_t:X^r\to X^r$ given by
$f_t(x,u)=(x,u+t)$ computed modulo identifications.
An invariant probability measure is given by $\mu^r=\mu\times\ell/\int r\,dm$
where $\ell$ is Lebesgue measure.

Next suppose that $f_t:\Omega\to\Omega$ is a semiflow preserving a
probability measure $\nu$.   Form the natural extension (or inverse limit)
$\hat\Omega$ consisting of curves $\{x(s)\}_{s\ge0}$ in $\Omega$
satisfying $f_t(x(s))=x(s-t)$ for all $s\ge t\ge 0$.
An invertible flow $\hat f_t:\hat\Omega\to\hat\Omega$ is given by
$\hat f_t\{x(s)\}_{s\ge0}=\{x(s+t)\}_{s\ge0}$.
The projection $\pi:\hat\Omega\to\Omega$, $\{x(s)\}_{s\ge0}\mapsto x(0)$
defines a semiconjugacy between $\hat f_t$ and $f_t$
and there is a unique $\hat f_t$-invariant measure $\nu$ on $\hat\Omega$
such that $\pi$ is measure preserving.

 \section{Mixing for geometric Lorenz flows}
 \label{sec-mix}
 
In this section, we prove Theorem~\ref{thm-geom}.
Namely, we show that if $f$ is an l.e.o.\ expanding map of Lorenz-type
and $r$ is a logarithmic roof function, then
the corresponding geometric Lorenz flow is mixing.
 By Remark~\ref{rmk-Ratner}, it is enough to prove that the
corresponding geometric Lorenz semiflow is weak mixing.

 Let $f:I\to I$ be a Lorenz-like expanding map.
As shown in Diaz-Ordaz~\cite{Ord04}, 
$f$ admits an {\em induced map} $F$ that is Gibbs-Markov.   More precisely 
there is an open interval $Y\subset I$ containing $0$ and
(modulo a set of Lebesgue measure zero) there is a partition
$\mathcal P=\{\omega\}$ of $Y$ consisting of intervals,
and a return function $R:Y\to\N$ constant on partition elements
such that the induced map $F(x)=f^{R(x)}(x)$ restricts
to a diffeomorphism $F|_\omega :\omega\to Y$ for each partition
element~$\omega$, and such that the following conditions are satisfied:
\begin{itemize}
\item[(a)] There exists \( \lambda > 1 \) such that
    $|DF_\omega|\geq \lambda$  for each $\omega$.
\item[(b)] For all \( \omega \),
  $\log g_\omega$ is H\"older (uniformly in $\omega$), where $g$ is the 
Jacobian of the inverse of $F|_\omega:\omega\to Y$.
\item[(c)] $R$ is Lebesgue integrable.
\end{itemize}
A standard Folklore Theorem in dynamics says that conditions~(a)--(c)
 implies the existence of a unique 
  ergodic invariant measure $\mu$ for $f:I\to I$ that is absolutely continuous 
  with respect to Lebesgue and whose support includes $Y$.
In addition, the construction of~\cite{Ord04} satisfies
\begin{itemize}
\item[(d)] $0\not\in \overline{f^k\omega}$ for $0\le k< R(\omega)$ for each 
$\omega$.
\end{itemize}

The following Liv\v{s}ic regularity result is due to 
Bruin~{\em et al.}~\cite{BruHolNic04}.

 \begin{lemma}  \label{lem-BHN}
Suppose that $f:I\to I$ admits an induced map $F:Y\to Y$ satisfying 
conditions~(a)--(d) above, with absolutely continuous ergodic measure $\mu$.   
Let $r:I\to\R^+$ be a H\"older roof function 
with a logarithmic singularity at $0$.
Let $\psi:I\to S^1$ be a $\mu$-measurable
function satisfying
\[
e^{ir}=(\psi\circ T)\psi^{-1},\enspace a.e.      
\]
Then $\psi$ has a version that is H\"older on $Y$.
\end{lemma}

\begin{proof}[Sketch proof]
We provide the details required for the reader to pass between the
formulation in~\cite{BruHolNic04} and the statement of the lemma.
By conditions~(a)--(c), the map $f:I\to I$ can be modelled by a Young 
tower~\cite{You99} with base $Y$.
The appropriate measure for the tower has a H\"older density, so
by condition~(d), Bruin~{\em et al.}~\cite{BruHolNic04}*{Theorem~2} applies to
guarantee that $\psi$ is H\"older on any fixed partition
element of $Y$ and hence on the whole of $Y$ (\cite{BruHolNic04}*{Remark~3}).
The tower metric in~\cite{BruHolNic04} coincides with the 
Euclidean metric when
restricted to $Y$, so $\psi$ is H\"older in the original metric.
\end{proof}

\begin{rmk}
 Related Liv\v{s}ic regularity results can be found in~\cite{Gou04c}.   
The arguments in~\cites{BruHolNic04,Gou04c} are different from each other
leading to distinct results, and the approach in~\cite{BruHolNic04} turns 
out to more convenient for our purposes.
%
\end{rmk}

Now let $r:I\to\R^+$ be a H\"older roof function with a logarithmic
singularity at $0$, and let $f_t$ be the corresponding geometric Lorenz
semiflow with ergodic measure $\mu^r$.

 \begin{lemma}  \label{lem-geom}
 Suppose that the geometric Lorenz semiflow is
 {\em not} weak mixing.   Then there exists a constant $a>0$
 and a measurable eigenfunction $\psi:X\to S^1$ continuous
 on $\bigcup_{k\ge0}f^k Y$,
 such that
 \begin{equation} \label{eq-mix}
 e^{iar}= (\psi\circ f) \psi^{-1} \enspace a.e.
 \end{equation}
 \end{lemma}

 \begin{proof}  It can be taken as the definition of weak mixing that
 $\phi\circ f_t=e^{iat}\phi$ has no measurable solutions 
$\phi:X^r\to S^1$ for $a>0$.
 Suppose that $\phi$ is such a measurable solution.
 It follows from Fubini
 that there exists $\epsilon>0$ with $r>\epsilon\;a.e.$ such that
 $\phi\circ f_t(x,\epsilon)=e^{iat}\phi(x,\epsilon)$ for almost every $x\in X$.
 Set $t=r(x)$ and $\psi(x)=\phi(x,\epsilon)$.   Since $f_{r(x)}(x,\epsilon)=(fx,\epsilon)$
 we obtain that $\psi\circ f=e^{iar}\psi$.
 Thus $\psi$ is a measurable solution to equation~\eqref{eq-mix}.
 By Lemma~\ref{lem-BHN}, there is a solution $\psi$ that
is H\"older continuous on $Y$.

 It is now straightforward to show that
$\psi$ is continuous on $\bigcup_{k\ge0}f^k Y$.
 Suppose that $z=f^k y$ where $y\in Y$.   Since $f(0)$ is undefined,
 it is certainly the case that $f^jy\neq0$ for $j=0,\ldots,k-1$.
 Hence we can choose an open set $U\subset Y$ containing $y$
 such that $0\not\in f^jU$ for $0\le j\le k-1$.
 In particular, there exists $\gamma\in(0,1)$ such that
 $\psi$ is $C^\gamma$ on $Y$ and at the same time
 $e^{iar_k}$ is $C^1$ on $U$.

 Let $z_i=f^k y_i$ where $y_i\in U$
 for $i=1,2$.    Applying equation~\eqref{eq-mix}
 in the form $\psi\circ f^k = e^{iar_k}\psi$, we obtain
 \[
 \psi(z_1)\psi(z_2)^{-1}=\psi(y_1)\psi(y_2)^{-1}e^{iar_k(y_1)}e^{-iar_k(y_2)}
 \]
 so by {\bf f(ii)},
 \[
 |\psi(z_1)\psi(z_2)^{-1}|\le D|y_1-y_2|^\gamma\le D(c\lambda^k)^{-\gamma} 
 |z_1-z_2|^{\gamma}.
 \]
 Hence $\psi$ is H\"older on $f^k U$ and so
 is certainly continuous at $z$ as required.
 \end{proof}

 \begin{proof}[Proof of Theorem~\ref{thm-geom}]
 Suppose that $f_t$ is not weak mixing.
 By Lemma~\ref{lem-geom},
 there exists $a>0$ and a measurable eigenfunction $\psi:I\to S^1$
 satisfying~\eqref{eq-mix}, such that $\psi$ is continuous on
 $\bigcup_{k\ge0}f^k Y$.
 Since $f$ is l.e.o., $\psi$ is continuous on $(-1,1)$.
In particular, $\psi$ is continuous at $f(\pm1)$.

Iterate equation~\eqref{eq-mix} to obtain
 \begin{align} \label{eq-2}
 e^{iar}e^{iar\circ f}=(\psi\circ f^2)\psi^{-1}.
 \end{align}
We evaluate this equation along a sequence
 $x_n>0$ with $x_n\to0$.  We have the limits
 $\psi(x_n)\to \psi(0)$, $\psi(f^2x_n)\to \psi(f(-1))$,
 and $r(fx_n)\to r(-1)$.
 Hence, the right-hand-side of equation~\eqref{eq-2} 
 and the second factor on the left-hand-side converge as $n\to\infty$.
 To obtain a contradiction, it suffices to choose $x_n\to0$ so that
 $e^{iar(x_n)}$ does not converge.    
 For $n$ sufficiently large, choose $b_n> r(\epsilon/2^n)$
 such that $e^{iab_n}=(-1)^n$.
 Since $r$ is continuous on $I-\{0\}$ 
 and $r(x)\to\infty$ as $x\to0$, there exists $x_n>0$
 such that $r(x_n)=b_n$.  As required, $x_n\to0$ and $e^{iar(x_n)}$
 does not converge.
 \end{proof}

\section{Extensions and future directions}
\label{sec-ext}

In this paper, we have focused on geometric Lorenz attractors satisfying
a locally eventually onto condition.
In this section, we discuss how to relax the l.e.o.\ condition, and
we describe how our results extend to larger classes of singular attractors.
We end by mentioning some future directions.

The l.e.o.\ assumption occurs in much of the 
literature since it is a verifiable property of the flow that guarantees
topological transitivity of $f$.   Moreover, it holds for the Lorenz
equations~\eqref{eq-lorenz} themselves~\cite{Tuc02}.
However, topological transitivity for $f:I\to I$ is not crucial in this paper.
We could simply work with the unique ergodic component whose support
includes the set $Y$ introduced in Section~\ref{sec-mix}.   
Passing to the natural extension of the suspension
over this component yields a geometric Lorenz attractor with an
SRB measure as before.

The proof of mixing for this attractor also relies to some extent on
the l.e.o.\ property, but it is clear that the proof goes through under a 
much weaker condition: namely that 
a preimage of a forward image of $1$ or $-1$ lies in $Y$.
(That is, there exist integers $k,\ell\ge0$ such that $f^\ell(1)\in f^kY$
or $f^\ell(-1)\in f^kY$.)   
This condition is open and dense within the class of Lorenz-like
expanding maps, and in principle, it can be verified by a finite computation.
 
Next, we note that the specific structure of Lorenz-like expanding maps
is not so crucial for the proof of mixing.   Properties~{\bf f(i)}--{\bf f(iii)}
could be relaxed, or altered completely, provided $f$
admits an induced map $F:Y\to Y$ satisfying conditions~(a)--(d).
The existence of a Gibbs-Markov induced map is a very general 
condition, see \cites{BruLuzStr03, AlvLuzPindim1}.  
 It clearly holds for uniformly expanding 
maps and, we have made use of the the recent work of Diaz-Ordaz~\cite{Ord04}
for Lorenz-like expanding maps.

Condition~{\bf f(ii)} in the definition of the Lorenz-like expanding
map corresponds to the eigenvalue condition
$\lambda_u>|\lambda_s|$ which is valid for the Lorenz attractor.
The resulting class of geometric Lorenz attractors are sometimes
called expanding.   There are also {\em contracting geometric Lorenz
attractors}~\cite{Rov93} that arise when $\lambda_u<|\lambda_s|$.
It seems likely that such attractors can be shown to be mixing using the 
ideas in this paper.

In situations where Lemma~\ref{lem-cone} fails, it may not be
possible to reduce to a one-dimensional expanding map.
However, it is plausible that the Poincar\'e map $P$ could be modelled
by an ``unquotiented'' Young tower as in~\cite{You98} to which the ideas
in this paper might still be applicable.

In this paper, we have established mixing.   An interesting open
question is to prove results on the speed of mixing.  In a different
direction, it would be interesting to derive statistical limit laws
such as the central limit theorem and invariance principles for
the Lorenz attractor.

\subsubsection*{Acknowledgements}
This research was supported in part by  EPSRC Grant GR/S11862/01.
IM is greatly indebted to the University of Houston for the use of e-mail,
given that pine is currently not supported on the University of Surrey network.

\end{document}